\documentclass{ifacconf}

\usepackage{graphicx}      
\usepackage{natbib}        

\usepackage[scientific-notation=false,binary-units=true]{siunitx}
\usepackage{multirow}

\usepackage{tabularx}
\newcolumntype{L}[1]{>{\raggedright\arraybackslash}p{#1}} 
\newcolumntype{C}[1]{>{\centering\arraybackslash}p{#1}} 

\usepackage{booktabs}
\usepackage{amsmath,amssymb,bm}
\usepackage{xcolor}

\usepackage[%
ngerman,
colorinlistoftodos,
textsize=scriptsize,	
]{todonotes}

\newcommand{\Autor}[1]{\citet{#1}} 
\newcommand{\Zitat}[1]{\citep{#1}} 
\newcommand{\Zitate}[1]{\citep{#1}} 

\newcommand{\Fig}[1]{Figure.~\ref{fig:#1}}
\newcommand{\fig}[1]{Fig.~\ref{fig:#1}}
\newcommand{\Tab}[1]{Tab.~\ref{tab:#1}}
\newcommand{\Eq}[1]{Eq.~\ref{eq:#1}}
\newcommand{\eg}{e.g.\ }
\newcommand{\ie}{i.e.\ }

\newcommand{\Matlab}{MATLAB}

\newcommand{\Kienzle}{\textsc{Kienzle}}
\newcommand{\KienzleMdl}{\Kienzle-model}
\newcommand{\Altintas}{\textsc{Altintas}}

\newcommand{\qml}{``}
\newcommand{\qmr}{''}

\newcommand{\grad}{\ensuremath{^\circ}}

\usepackage[%
nolist,
]{acronym}

\definecolor{wzl-blue}			{RGB}{0	, 84,199}
\definecolor{wzl-LightBlue}		{RGB}{199,221,242}
\definecolor{wzl-DarkGray}		{RGB}{128,128,128}
\definecolor{wzl-LightGray}		{RGB}{191,191,191}
\definecolor{wzl-red}			{RGB}{204,  7, 30}
\definecolor{wzl-yellow}		{RGB}{246,168,  0}
\definecolor{wzl-green}			{RGB}{  0,153, 51}

\usepackage{tikz,pgfplots}
\usetikzlibrary{spy}
\usepgfplotslibrary{fillbetween} 

\usetikzlibrary{external}
\tikzexternaldisable

\pgfplotsset{compat=newest} 
\makeatletter

\pgfplotsset{
	boxplot prepared from table/.code={
		\def\tikz@plot@handler{\pgfplotsplothandlerboxplotprepared}%
		\pgfplotsset{
			/pgfplots/boxplot prepared from table/.cd,
			#1,
		}
	},
	/pgfplots/boxplot prepared from table/.cd,
	table/.code={\pgfplotstablecopy{#1}\to\boxplot@datatable},
	row/.initial=0,
	make style readable from table/.style={
		#1/.code={
			\pgfplotstablegetelem{\pgfkeysvalueof{/pgfplots/boxplot prepared from table/row}}{##1}\of\boxplot@datatable
			\pgfplotsset{boxplot/#1/.expand once={\pgfplotsretval}}
		}
	},
	make style readable from table=lower whisker,
	make style readable from table=upper whisker,
	make style readable from table=lower quartile,
	make style readable from table=upper quartile,
	make style readable from table=median,
	make style readable from table=lower notch,
	make style readable from table=upper notch,
	make style readable from table=draw position,
}
\makeatother

\pgfplotstableset{col sep=comma}

\newcommand{\SigFromTableSingle}[3]{
	{
		\addplot[mark=none,color=#3] table [x=x, y=#2, col sep=comma] {#1};		
	};
}

\newcommand{\SigFromTable}[1]{
	\pgfplotstableread[col sep = comma]{#1}\mydata
	\pgfplotstablegetcolsof{\mydata}
	\pgfmathtruncatemacro\TotalCols{\pgfplotsretval-1}
	\pgfplotsinvokeforeach{1,...,\TotalCols}
	{
		\addplot+[mark=none] 				table[x=x, y index=##1 , col sep=comma] {\mydata};
		\addlegendentry{##1};
	};
}

\newcommand{\SigLbUbFromTableB}[2]{
	{synctex=1
		\addplot[mark=none,color=#2,draw=none,name path=lb] 	table [x=x, y=lb,  col sep=comma] {#1};
		\addplot[mark=none,color=#2,draw=none,name path=ub] 	table [x=x, y=ub,  col sep=comma] {#1};
		\addplot[
		thick,
		color=#2,
		fill=#2, 
		fill opacity=0.2
		]	fill between[of=lb and ub];
		\addplot[mark=none,thick,color=#2] 				table [x=x, y=sig, col sep=comma] {#1};
	};
}

\DeclareRobustCommand\sampleline[1]{%
	\tikz\draw[#1] (0,0) (0,\the\dimexpr\fontdimen22\textfont2\relax)
	-- (2em,\the\dimexpr\fontdimen22\textfont2\relax);%
}

\pgfplotstableset{col sep=comma}
\pgfplotsset{major grid style={dotted,wzl-LightGray}}
\pgfplotsset{minor grid style={dotted,wzl-LightGray}}

\def\NumSlices{23}
\def\NumRev{10} 
\def\NumSimKnzl{1000}

\begin{document}
	\begin{frontmatter}
		
		\title{Identifying trending coefficients with an ensemble Kalman filter} 
		
		
		
		\author[WZL] {M. Schwenzer} 
		\author[igpm]{G. Visconti}
		\author[IRT] {M. Ay}
		\author[WZL] {T. Bergs}
		\author[igpm]{M. Herty}
		\author[IRT] {D. Abel}
				
		\address[WZL]{Laboratory for Machine Tools and Production Engineering (WZL), RWTH Aachen University, Campus-Boulevard 30, 52074 Aachen, Germany (e-mail: m.schwenzer@wzl.rwth-aachen.de).}
		\address[igpm]{Institute of Geometry and Applied Mathematics 
			(IGPM), RWTH Aachen University,  Templergraben 55 , 52056  Aachen, Germany}
		\address[IRT]{Institute of Automatic Control, RWTH Aachen University,  Campus-Boulevard 30, 52074 Aachen, Germany}

		\begin{abstract}                
			This paper extends the \ac{EnKF} for inverse problems to identify trending model coefficients. This is done by repeatedly inflating the ensemble while maintaining the mean of the particles.
			As a benchmark serves a classic \ac{EnKF} and a \ac{RLS} on the example of identifying a force model in milling, which changes due to the progression of tool wear.
			For a proper comparison, the true values are simulated and augmented with white Gaussian noise. The results demonstrate the feasibility of the approach for dynamic identification while still achieving good accuracy in the static case. Further, the inflated \ac{EnKF} shows a remarkably insensitivity on the starting set but a less smooth convergence compared to the classic \ac{EnKF}.			
			
		\end{abstract}
		

\end{frontmatter}

\begin{acronym}[EnKF]		
	\acro{MPC}[MPC]{model-based predictive control}	
	
	\acro{KF}[KF]{Kalman filter}
	\acro{EKF}[EKF]{extended Kalman filter}
	\acro{EnKF}[EnKF]{ensemble Kalman filter}
	\acro{RLS}[RLS]{recursive least squares}
	
	\acro{rms}[rms]{root mean square}
	\acro{snr}[S/N]{signal-to-noise ratio}
	\acro{KKT}{Karush-Kuhn-Tucker}
\end{acronym}

\section{ Introduction } 
%
In milling, a rotating cutting movement is overlaid with a translatory feed movement resulting in an cyclically intermittent cutting process. 
It is a flexible and highly dynamic manufacturing process for free-form surfaces. This leads to a continuously varying thickness of the removed chip and with it a varying force. 

The force is the most important parameter to analyze and evaluate the cutting process. Ever since there were machine tools, researchers aimed at describing the force through models in order to better understand and design the process. 
Nowadays, those models are also used for advanced control of the milling process, be it adaptive control \Zitat{altintas_integration_2017} or even \ac{MPC} \Zitat{stemmler_model_2017}. 
The semi-empirical models must be calibrated for every tool-workpiece material combination.
They represent a certain tool state and change as the tool wears. 
Historically, the identification was conducted off-line but new approaches have paved the way for an on-line identification in the manufacturing process.\newline
This work discusses the problem of identifying changing coefficients. For this a constrained \ac{EnKF} is repeatedly inflated and benchmarked against a \ac{RLS}. The overall objective is to present a method for continuous identification of time-variant models.

\section{ State of the art } 

\subsection{ Force model }
So-called mechanistic models relate the force to the cross-section of the undeformed chip.
The most popular examples are the exponential force model according to \Autor{kienzle_bestimmung_1952} and the linear approach according to \Autor{altintas_general_1996}:
\begin{align}
	F_{i,Kienzle}  &= k_{i} \, b \, h^{1-m_i} \hspace{1cm} &i\in t,r, \label{eq:KienzleForceModel}\\
	F_{i,Altintas} &= b \, \left(K_{i,e} + K_{i,c}\, h\right) &i\in t,r. \label{eq:AltintasForceModel}
\end{align}
The parameters $k_{i}$ and $m_i$, or $K_{i,e}$ and $K_{i,c}$ respectively, represent material-specific coefficients. The indices $i = t,r$ indicate the tangential and the radial component of the force vector. The undeformed chip thickness $h$ and the undeformed chip width $b$ form the cross-section of the chip.

Often, the admissibility range of linear models is extended by assuming the coefficients to be themselves again a function of the undeformed chip thickness $h$, \eg as a linear \Zitat{grossi_accurate_2017}, an exponential \Zitate{wan_new_2007,wan_new_2009,campatelli_prediction_2012,zhang_non-contact_2018}, or a polynomial relation \Zitate{wei_prediction_2018,wang_comparison_2018}. 
This converts the linear model with varying coefficients $K = f(h)$ into a non-linear -- often exponential -- model with constant coefficients \Zitate{wan_new_2007,yao_chatter_2013}.

\subsection{ Model identification }
The great majority of the work on how to identify the coefficients of mechanistic force models has been done for linear models, namely the force model of \Altintas{} (\Eq{AltintasForceModel}). 
Few works identify a non-linear force model  \Zitate{jayaram_estimation_2001,dotcheva_evaluation_2008,wang_identification_2013,adem_identification_2015,zhang_new_2017,zhang_non-contact_2018} or the \KienzleMdl{} explicitly \Zitate{shin_new_1997,perez_enhanced_2013,schwenzer_comparative_2018}.\

There exist two ways to identify such mechanistic force models in milling: the method of average force and the method of instantaneous undeformed chip thickness.
The first is inspired from turning where the chip geometry does not vary resulting in a static cutting force. In milling this is approximated by averaging the dynamic force signal over a revolution. The method requires several dedicated experiments and is not on-line capable. \newline
The method of instantaneous undeformed chip thickness is essentially a curve fit between measurements and a simulated force signal. 
They have been dominated by global optimization, such as evolutionary algorithms \Zitate{grossi_accurate_2017,chen_precise_2018}, or particle swarm algorithms \Zitat{zhang_new_2017} for both presented models (\Eq{KienzleForceModel}, \Eq{AltintasForceModel}).
Nevertheless, local optimization algorithms have a significant advantage in computation with no loss in accuracy \Zitat{freiburg_determination_2015,schwenzer_comparative_2018}.

\Autor{adem_identification_2015} compared both models and identification methods. They concluded that a non-linear force model is generally more accurate than a linear model and that the optimization-based curve fit results in more accurate models than the average forces approach. 
\Autor{gonzalo_method_2010} focused on the latter comparing identification methods. They define a \qml{}true\qmr{} reference by identifying the coefficients in turning. Though the improvement is small, they argue that the method of instantaneous undeformed chip thickness has a better physical credibility due to the correspondence to the turning coefficients.

First studies on a continuous formulation of the method of instantaneous undeformed chip thickness propose an \ac{EnKF} as a non-linear estimator to identify the \KienzleMdl{} \Zitate{schwenzer_ensemble_2019,schwenzer_continuous_2019}. The studies revealed the extraordinary insensitivity of the \ac{EnKF} to measurement noise. They used an unscaled version of the filter as a trade-off for fast convergence against stability. 


\section{ Approach } 
Progressive tool wear changes the model coefficients gradually. Therefore, we examine the case of trending coefficients as a special case of changing coefficients. 
For this, we suggest three different approaches for a continuous identification of a time-varying model:
\begin{itemize}
	\item \ac{RLS},
	\item classic \ac{EnKF},
	\item \ac{EnKF} with repeated ensemble inflation (\ac{EnKF}$^\star$).
\end{itemize}
In order to meet generality, we use the non-linear \KienzleMdl{}, which is in fact not observable. Neither the \ac{RLS}, nor the \ac{EnKF} require observability of the system as they do not demand uniqueness of the solution of the inverse problem.

\subsection{ Recursive least square } 
The \ac{RLS} algorithm is a least square fit taking the estimate of the previous time instance into account. 
It works as similar to an exponential smoothing with a decreased weighting of the previous information.

Assume an arbitrary measurement system 
\begin{align}
	\hat{\bm{y}}_k & = f(\hat{\bm{x}}_k) \approx \bm{M}_k^\intercal \, \hat{\bm{x}}_k \text{,}
\end{align}
with a measurement $\hat{\bm{y}}_k$ that has non-linear relationship to the estimated state vector $\hat{\bm{x}}_k$ for time instance $k$ \Zitat{strejc_least_1979}. For a proper state-space representation, the relationship is linearized through the measurement-matrix $\bm{M}_k^\intercal$.
The model-matrix weights the error between the measurement and the prediction through the gain $\bm{G}_k$
\begin{align}
	\hat{\bm{x}}_k &= \hat{\bm{x}}_{k-1} + \bm{G}_k \, \left(\bm{z}_k - f(\hat{\bm{x}}_k)\right),\\
	\bm{G}_k  &= \frac{ \bm{P}_{k-1} \, \bm{M}_k }{ \rho + \bm{M}_k^\intercal \, \bm{P}_{k-1} \, \bm{M}_k },\label{eq:KalmanGain}\\
	\bm{P}_k  &= \frac{1}{\rho} \, \left[ \bm{I} - \bm{G}_k \, \bm{M}_k^\intercal \right]\, \bm{P}_{k-1}.\label{eq:Covariance}
\end{align}
The recursive approximation of the covariance matrix $\bm{P}_k$ limits the computational complexity. In the case on hand, the state vector becomes $\hat{\bm{x}}_k = \left[k_i \quad m_i\right]^\intercal$.
The forgetting factor $0 < \rho \leq 1$ (here $\rho = 0.98$) weights the samples.
The initial value of the covariance matrix is set to $\bm{P}_0 =  10^5 \, \bm{I}$, with the identity matrix $\bm{I}$. Applying this to

\subsection{ Ensemble Kalman filter } 
Instead of integrating a single state vector forward in time, the \ac{EnKF} propagates an ensemble of state vectors and takes its mean as the best-guess \Zitat{evensen_sequential_1994}. 
The backbone of the \ac{EnKF} is a Markov chain  Monte Carlo simulation of the evolution in time of individual state vectors, which approximate the true probability density of the states \Zitat{evensen_ensemble_2003}. The \ac{EnKF} is a special case of a particle filter without re-sampling and with Gaussian distribution for the measurement likelihood.

Assuming an ensemble matrix 
\begin{equation}
	\bm{X}_{k|k-1} = (\bm{x}_{k|k-1}^1, \dots, \bm{x}_{k|k-1}^J),
\end{equation}
with $J$ individual state vectors $\bm{x}_{k|k-1}^j,\ j\in \{1,\dots,J\}$ at time step $k$ based on the information from time step $k-1$.  In the case on hand, the state vector is
\begin{equation}
	\begin{tabular}{l c r r r r r r}
		$\bm{x}$ 	 &= [ & $F_t$    & $F_r$ 	&$k_t$& $k_r$&$m_t$&$m_r ]^\intercal$,
	\end{tabular}
\end{equation}
consisting out of the measurable forces $F_i$ and the non-observable parameters $k_i$, $m_i$. The \ac{RLS} does not require the measurements to be a member of the state vector.

The \ac{EnKF} is a truly non-linear estimator, using the model function $f(\bm{x},\bm{u})$ to propagate every member of the ensemble forward in time based on the current input $\bm{u}_k$:
\begin{align}
	\hat{\bm{X}}_{k|k-1} &= f\left( \bm{X}_{k-1|k-1},\bm{u}_k \right),\\
	\bm{P}_{k|k-1} &= \frac{1}{J}~ (\hat{\bm{X}}_{k|k-1} - \overline{\bm{X}}_{k|k-1})\, (\hat{\bm{X}}_{k|k-1} - \overline{\bm{X}}_{k|k-1})^\intercal. \label{eq:EnKF_Pkk-1}
\end{align}
The error covariance matrix is approximated by the empirical covariance matrix of the ensemble $\bm{P}$. It represents the spread to the mean of the ensemble $\overline{\bm{X}}$. 
Recent studies \Zitate{schwenzer_ensemble_2019,schwenzer_continuous_2019} resign from scaling the covariance $\bm{P}$, which favors the convergence speed but comes at the cost of a higher instability. The unscaled version of the filter looses the feature of mathematical well-posedness, that is why it is not further considered.


The predicted state vectors are updated through the measurements weighted by the Kalman gain $\bm{G}_k$:
\begin{align} 
	\bm{G}_k 	&= \bm{P}_{k|k-1} \, \bm{H}_k^\intercal ~ \left[ \bm{H}_k \, \bm{P}_{k|k-1} \, \bm{H}_k^\intercal + \overline{\bm{R}}_k \right] ^{-1}, \label{eq:EnKF:Gain} \\  
	\bm{x}_{k|k}^j &= \hat{\bm{x}}_{k|k-1}^j + \bm{G}_k \left[ \bm{z}_k^j -  %
	\bm{H}_k \,\hat{\bm{x}}_{k|k-1}^j%
	\right]
\end{align}
where $\bm{H}_k$ is the linear observation operator and the single measurement vector $\bm{y}_k$  is augmented with artificial zero-mean Gaussian noise  $\bm{\varepsilon}\sim \mathcal{N}(\bm{0},\,\bm{\Gamma})$:
\begin{align}
	\bm{z}_k^j 	&= \bm{y}_k + \bm{\varepsilon}_k^j,\qquad j \in \{1,\dots,J\} \label{eq:EnKF_MeasurementNoise_z}\\
	\bm{Z}_k 	&= \left(\bm{z}_k^1,\dots, \bm{z}_k^n\right)^\intercal,\quad \text{and}\\
	\overline{\bm{R}}_k 	&= \frac{1}{J}\, \left(\bm{\varepsilon}_k \, \bm{\varepsilon}_k^\intercal\right). \label{eq:EnKF:MeasurementNoise_zR}
\end{align}
The update of the empirical covariance matrix $\bm{P}_{k|k}$ is calculated as before, but now uses the information of the current time step $k$
\begin{align}
	\bm{P}_{k|k} &= \frac{1}{J}\, (\bm{X}_{k|k} - \overline{\bm{X}}_{k|k})\, (\bm{X}_{k|k} - \overline{\bm{X}}_{k|k})^\intercal.\label{eq:EnKF_Pkk}
\end{align}
Note that the Kalman gain $\bm{G}_k$, \Eq{EnKF:Gain}, and the covariances $\bm{P}_{k|k-1}$, \Eq{EnKF_Pkk-1}, and $\bm{P}_{k|k}$, \Eq{EnKF_Pkk}, are the same for all ensemble members and only need to be calculated once. 

In the case of application to inverse problems, the \ac{EnKF} can be written in the following formulation
\begin{align}
	\bm{X}_{k+1|k+1} =& \bm{X}_{k|k} + \bm{C}_{k|k} \, \left(\bm{D}_{k|k} + \bm{\Gamma}\right)^{-1} \, \left( \bm{Z}_{k+1}-\hat{\bm{Y}}_{k} \right),\label{eq:SdT:EnKF_math_update}\\
	\bm{C}_{k|k} =& \frac{1}{J} (\hat{\bm{Y}}_k - \overline{\hat{\bm{Y}}}_k)\,(\hat{\bm{Y}}_k - \overline{\hat{\bm{Y}}}_k)^\intercal,\\
	\bm{D}_{k|k} =& \frac{1}{J} (\bm{X}_{k|k} - \overline{\bm{X}}_{k|k})\,( %
	\hat{\bm{Y}}_k%
	- \overline{\hat{\bm{Y}}}_k	)^\intercal,
\end{align}
which is often used in the field of mathematics in order to analyse the general behavior of the filter and to provide a comprehensive theoretical background.
Whereas often the measurements $\bm{y}_k$ are not augmented: $\bm{\varepsilon}_k = \bm{0} \Rightarrow \bm{Z}_k = \bm{y}_k$. This facilitates the analysis of the stability and well-posedness of the \ac{EnKF}.
However, perturbing the measurements creates individual observations for every member and prevents particles from synchronizing and collapsing completely to a potentially non-optimal solution \Zitat{kelly_well-posedness_2014}.

Although the filter converges within the subspace spanned by the initial ensemble -- \ie the mean value of the ensemble converges to a value within this subset -- it is not guaranteed that all ensemble members stay within this subset at all times. This makes it fragile if the model used is not defined for all states or exhibits singularities.\newline
Enforcing box constraints, hoping that the particular member does not get saturated but will (eventually) move back into the subspace some time, is reasonable if the initial ensemble is chosen properly. That is if the solution lays within the limits, so that the whole ensemble is drawn towards its mean by the common covariance. Further, if the ensemble is large enough, it is unlikely that many members are affected by artificial box constraints and; thus, distort the expected behavior of the filter.
In fact, \Autor{chada_incorporation_2019} showed that the \ac{EnKF} can impose box-constraints by projecting states that violate a constraint back onto the boundary (this is the projected Newton method). This might affect the convergence because it changes  the direction and the step-size of the particle. Considering the \ac{EnKF} as a sequential optimization, it becomes obvious that the optimization gets affected. Nevertheless, it is still ensured that the ensembles collapse to their mean and that the method converges to a \ac{KKT} point of the optimization problem.

In the case of identifying mechanistic force models in milling, the box-constraints for the state vector $\bm{x}$ were set to
\begin{equation}
	\begin{tabular}{l c r r r r r r}
		$\bm{x}_{lb}$ &= [ &$-\infty$ &$-\infty$ & 500 & 100  &0.1  &0.1 $]^\intercal$,\\
		$\bm{x}_{ub}$ &= [ &$ \infty$ &$ \infty$ &3500 &2100  &  1  &  1 $]^\intercal$.
	\end{tabular}
\end{equation}

\subsection{ Ensemble Kalman filter with repeated ensemble inflation } 
The problem of identifying time-variant models could be approached by restarting the \ac{EnKF} if the model accuracy runs out of bound. This would neglect all prior gained information and reset the filter.
In contrast to restarting the \ac{EnKF} over and over again, the ensemble inflation should maintain the current mean $\overline{\bm{X}}_k^\star \approx \overline{\bm{X}}_k$ but be brought back to the initial spread, \ie covariance $\bm{P}_{0}$.

On a regular basis, the ensemble $\hat{\bm{X}}_k$ is substituted by an ensemble $\hat{\bm{X}}_k^\star$ with the covariance 
\begin{equation}
	\bm{P}_{k|k}^\star = \frac{1}{\lambda}\, \bm{P}_{0} \label{eq:P_inflt}
\end{equation}
In contrast to this, variance inflation scales the covariance ahead of the analysis, \Eq{EnKF:Gain} or \Eq{SdT:EnKF_math_update} respectively. It keeps the particles from collapsing too fast, which is of importance at problems of high dimensions (large number of states) and a comparatively small ensemble size.\newline
Ensemble inflation replaces the complete ensemble in discrete steps after analysis.

The idea of the mean field theory is that many small random subsets of particles on average represent the behavior of the whole particle cloud \Zitat{herty_kinetic_2019}. 
Following the mean field theory, we inflate only a small subset of the ensemble set $\bm{\mathcal{X}} \subset \{\bm{x}^j\}_{j=1}^J$, $|\bm{\mathcal{X}}| = M \ll J$, with regard to its mean. This should smooth the inflation or the immediate effect of the inflation respectively.
The ensemble size was set to $J = 100$ and the subset $M$ to \SI{10}{\percent}.

\subsection{ Simulation for validation } 
In general, three distinct cases must be considered for identification:
\begin{itemize}
	\item static coefficients,
	\item ascending coefficients, and
	\item alternating coefficients.
\end{itemize}
For validation, this paper simulates \NumRev{} revolutions of a straight milling operation with a flat end-mill emulating a machining of stainless steel (X5CrNi18-10) with a two-fluted solid carbide tool, \Tab{ProcessParameters}. Those are the same conditions as in \Zitate{schwenzer_ensemble_2019,schwenzer_continuous_2019} for reasons of consistency and a better comparison. 

\begin{table}
	\centering
	\caption{Tool and process parameters} \label{tab:ProcessParameters}
	\begin{tabularx}{8.855cm}{L{1.7cm} C{.3cm} C{1cm} || L{1.7cm} C{.3cm} C{1.2cm}}
		\toprule
		\multicolumn{3}{C{3cm}||}{Tool geometry}			& \multicolumn{3}{C{3cm}}{Process parameter}	\\ \midrule
		Diameter		& $D$		& \SI{10}{\milli\meter} & Feed				& $f$		& \SI{0.1}{\milli\meter}\\
		Number of teeth & $N_z$		& 2 					& Cutting velocity	& $v_c$ 	& \SI{2.44}{\meter\per\second}	\\
		Helix angle 	& $\beta$ 	& \SI{45}{\grad} 		& Depth of cut		& $a_p$		& \SI{2}{\milli\meter}	\\
		Rake angle 		& $\gamma$	& \SI{20}{\grad}		& Width of cut		& $a_e$		& \SI{3}{\milli\meter}	\\
		\bottomrule
	\end{tabularx}
\end{table}

Since mechanistic force models are defined for straight cutting edges, the helix angle of the tool is approximated as a spiral staircase, in this case with \NumSlices{} disk elements. 
The force measurements were simulated assuming the coefficients given in \Tab{KienzleCoefs_theory_X5}, which were taken from one of the most extensive work on coefficients of the \KienzleMdl{} \Zitat{konig_spezifische_1982}. Using a sample frequency of $f_s  = \SI{10}{\kilo\hertz}$, a revolution translates  to approximately $234$ samples. Artificial Gaussian noise with \ac{snr} of 15 was added to the simulated force signals in order to mimic severe measurement noise. \fig{F_sim} illustrates the different cases and the added noise.
Simulating the force signals allows to use the coordinate system of the cutting edge, in which the model is defined, and to evaluate if the identification tends towards the true coefficients.

\begin{table}	
	\centering
	\caption{Force coefficients for the Kienzle force model for X5CrNi18-10}
	\label{tab:KienzleCoefs_theory_X5}
	\begin{tabular}{p{3.8cm} | c  c  c  c}
		\toprule
		& $k_{t}$	& $m_t$	& $k_{r}$	& $m_{r}$ \\ \midrule
		upper bound	& 1800		& 0.6	& 1200			& 0.3 	\\
		lower bound	& 800		& 0.05		& 600		& 0.01  	\\ 
		X5CrNi18-10 \cite[s.~95]{konig_spezifische_1982} & 1700 	& 0.18 & 350 & 0.55\\
		\bottomrule
	\end{tabular}
\end{table}			

\def\SzW{\columnwidth}
\def\SzH{3.5cm}

\def\xmin{0}
\def\xmax{2300}
\def\ymin{-75}
\def\ymax{700}

\setlength{\tabcolsep}{0pt} 
\tikzexternalenable

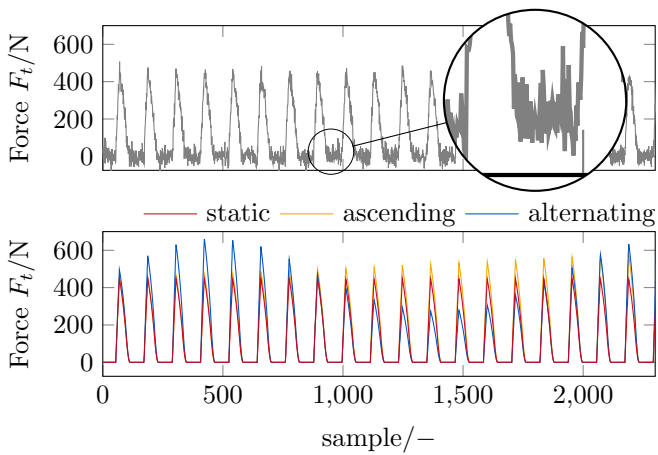
\begin{figure}[htp]
	\begin{tabular}{ c }
		\tikzsetnextfilename{F_sim_noise} 
		\begin{tikzpicture}[spy using outlines={circle, magnification=4, connect spies}]
			\begin{axis}[
							xmin=\xmin,		xmax=\xmax,
							ymin=\ymin,		ymax=\ymax,
							ylabel={Force $F_t / \si{\newton}$},
							xticklabels=\empty,
							width=\SzW,		height=\SzH,
							cycle list={{wzl-DarkGray},{wzl-blue},{wzl-red}},
							]
				\SigFromTable{F_wgn.csv}
				\legend{}; 
				
				\coordinate (spypoint) at (axis cs:950,25);
				\coordinate (magnifyglass) at (axis cs:1800,300);
			\end{axis}
			
			\spy [black, size=2.4cm] on (spypoint) in node[fill=white] at (magnifyglass);
		\end{tikzpicture} \\
	\multicolumn{1}{r}{%
			\sampleline{wzl-red} static \sampleline{wzl-yellow} ascending \sampleline{wzl-blue} alternating%
		}\\
		\tikzsetnextfilename{F_sim_ideal} 
		\begin{tikzpicture}
			\begin{axis}[
						xmin=\xmin,		xmax=\xmax,
						ymin=\ymin,		ymax=\ymax,
						ylabel={Force $F_t / \si{\newton}$},
						xlabel={$\text{sample} / -$},
						width=\SzW,		height=\SzH,
						cycle list={{wzl-yellow},{wzl-blue},{wzl-red}},
						]
				\SigFromTable{F_act.csv}
				\legend{}; 
			\end{axis}
		\end{tikzpicture}\\ 

	\end{tabular}

\caption{
		Simulated force signals: added white Gaussian noise with \ac{snr} = 15 for the static case (top); different cases without Gaussian noise (bottom)
		}
	\label{fig:F_sim}
\end{figure}

\tikzexternaldisable
\setlength{\tabcolsep}{6pt} 

The \ac{EnKF} embraces randomness as it depends on the choice of the initial ensemble and of the random perturbations at every step. In order to perform a sensitivity analysis, the \ac{EnKF} is restarted with $N_{init} = \NumSimKnzl{}$ different initial ensembles for every simulation case. The initial ensembles were generated as a uniform random distribution within the boundaries given in \Tab{KienzleCoefs_theory_X5}. They were generated once beforehand and the pseudo-random generator is initialized to the same value for every simulation case. 
For the \ac{RLS}, the mean of the initial ensemble served as the starting value.

When the undeformed chip thickness is tool small, cutting turns into ploughing and models loose their validity. Therefore, the identification was set to transit when the undeformed chip thickness (the sum along the disk elements) was smaller than a threshold $h_{TH} = \SI{0.01}{\milli\meter}$. Those corrected samples are denoted as samples$^\prime$.

For the sake of conciseness and without loss of generality, the discussion is limited to the tangential component of the force and the coefficients $i = t$. 
All calculations were performed with the software \Matlab{} R2018b from The MathWorks on an AMD Ryzen7-2700 (\SI{3}{\giga\hertz}) computer running Windows 10.

\section{ Results and discussion } 
From a grid-search we learned that large steps and a small inflation (large $\lambda$) are better for slowly changing model coefficients; while small steps and a large inflation is required for dynamic and drastic changes, \Tab{RMS}. Though this finding is intuitive, it is difficult to choose the optimal parameters. A good compromise were small steps and a small inflation because the increase in accuracy for quasi-static problems is much smaller than the loss in accuracy for a dynamic problem. Therefore, the step-size, in which the ensemble is inflated, was set to 50 samples$^\prime$ and the inflation factor to $\lambda = 10$.
\begin{table}
	\centering
	\caption{ Root-mean-square of the error, \Eq{eF}, for different inflation steps and fractions } \label{tab:RMS}
	\begin{tabular}{l | c | c || c | c | c }
		\toprule		
		Method 						& step 		& $\lambda$ & static& ascending & alternating 	\\
		\midrule
		\ac{EnKF}$^\star$			&	50		&	1.5		& 	7.6	&	8.4		& 8.5	\\
		\ac{EnKF}$^\star$			&	50		&	1		&	7.8	&	8.6		& 8.6	\\
		\ac{EnKF}$^\star$			&	50		&	10		&	6.4	&	7.2		& 8.9	\\
		\ac{EnKF}$^\star$			&	50		&	2		&	7.5	&	8.3		& 8.5	\\
		\ac{EnKF}$^\star$			&	50		&	5		&	6.9	&	7.7		& 8.4	\\
		\ac{EnKF}$^\star$			&	100		&	1.5		&	6.5	&	7.2		& 11.9	\\
		\ac{EnKF}$^\star$			&	100		&	1		&	6.6	&	7.3		& 11.5	\\
		\ac{EnKF}$^\star$			&	100		&	10		&	5.6	&	6.6		& 14.1	\\
		\ac{EnKF}$^\star$			&	100		&	2		&	6.4	&	7.1		& 12.2	\\
		\ac{EnKF}$^\star$			&	100		&	5		&	6.0	&	6.8		& 13.0	\\
		\ac{EnKF}$^\star$			&	200		&	1.5		&	4.6	&	7.3		& 20.2	\\
		\ac{EnKF}$^\star$			&	200		&	1		&	4.7	&	7.3		& 19.6	\\
		\ac{EnKF}$^\star$			&	200		&	10		&	4.2	&	7.3		& 23.6	\\
		\ac{EnKF}$^\star$			&	200		&	2		&	4.6	&	7.3		& 20.6	\\
		\ac{EnKF}$^\star$			&	200		&	5		&	4.4	&	7.2		& 22.2	\\
		\ac{EnKF}					&$\infty$ 	&	-		&	3.8	&	23.2	& 50.7	\\
		\bottomrule
	\end{tabular}
\end{table} 

\Fig{eFt_evol} illustrates the evolution of the error for the different identification methods in all three cases. The transparent tubes indicate the variance within the \NumSimKnzl{} simulations with different initial ensembles ($\overline{\Delta\bm{F}_t} \pm 2\,\sigma(\Delta\bm{F}_t)$). That are \SI{95.45}{\percent} if this Monte Carlo simulation is normally distributed. The error is defined as the difference in the force between the identified and the ideal coefficients
\begin{equation}
	\Delta F_t = F_{t,sim} - F_{t,ideal} = \hat{k}_t \, b \, h^{1-\hat{m}_t}  - k_t \, b \, h^{1-m_t}. \label{eq:eF}
\end{equation}
One can see that the \ac{RLS} (yellow lines) hardly depends on the initial value. Therefrom, it is remarkable how drastic the variation of the mean became; in particular when identifying static coefficients. The error even ran out of the box and broke down in a singularity as no box-constraints were imposed. In general, the \ac{RLS} exhibited a high oscillation in the mean error (thick lines) with no sign of convergence. Previous studies suggested that it is particularly prone to measurement noise \Zitat{schwenzer_continuous_2019}.\newline
The dynamic cases (\qml{}ascending\qmr{} and \qml{}alternating\qmr{}) revealed that the classic \ac{EnKF} (blue lines) was not sufficient but only worked ideal in the case of static model coefficients. However, inflating the ensemble repeatedly (red) decreased the accuracy in the static case but was the only option to achieve accurate results for identifying trending coefficients. The samples where the ensemble was inflated are marked by the dotted grid-lines.

\def\SzW{8cm}
\def\SzH{3cm}

\def\xmin{-10}
\def\xmax{1150}
\def\ymin{-60}
\def\ymax{60}
\def\xTick{0,200,400,600,800,1000,1200}
\def\NumMinorXtick{4}
\def\ColorA{wzl-yellow}
\def\ColorB{wzl-blue}
\def\ColorC{wzl-red}
\def\ColorD{black}

\setlength{\tabcolsep}{0pt} 
\tikzexternalenable

\begin{figure}[htp]
	\begin{tabular}{ l  c } 
		& \multicolumn{1}{r}{
			\sampleline{\ColorA} \ac{RLS} \sampleline{\ColorB} \ac{EnKF} \sampleline{\ColorC} \ac{EnKF}$^\star$
		}
		\\
		\tikzsetnextfilename{eFt_evol_sta} 
		\raisebox{0.7cm}{	\rotatebox[origin=c]{90}{ static }		}
		&
		\begin{tikzpicture}
			\begin{axis}[
						ymin=\ymin,		ymax=\ymax,
						ylabel={$\Delta F_t / \si{\newton}$}, 
						xmin=\xmin,		xmax=\xmax,
						xtick=\xTick,	xticklabels=\empty,
						width=\SzW,		height=\SzH,
						xminorgrids,	minor x tick num=\NumMinorXtick,
						xmajorgrids,
						]
				\SigLbUbFromTableB{eF_RLS_sta.csv}{\ColorA}

				\SigLbUbFromTableB{eF_xEnKF_sta.csv}{\ColorC}
				\SigLbUbFromTableB{eF_EnKF_sta.csv}{\ColorB}

			\end{axis}
		\end{tikzpicture}
		\\
		\hline 
		\raisebox{0.7cm}{	\rotatebox[origin=c]{90}{ ascending }		}
		&
		\tikzsetnextfilename{eFt_evol_asc} 
		\begin{tikzpicture}
			\begin{axis}[
						ymin=\ymin,		ymax=\ymax,
						ylabel={$\Delta F_t / \si{\newton}$}, 
						xmin=\xmin,		xmax=\xmax,
						xtick=\xTick,	xticklabels=\empty,
						width=\SzW,		height=\SzH,
						xminorgrids,	minor x tick num=\NumMinorXtick,
						xmajorgrids,
						]
				\SigLbUbFromTableB{eF_RLS_asc.csv}{\ColorA}
				\SigLbUbFromTableB{eF_EnKF_asc.csv}{\ColorB}
				\SigLbUbFromTableB{eF_xEnKF_asc.csv}{\ColorC}

			\end{axis}
		\end{tikzpicture}
		\\
		\hline 
		
		\raisebox{1.4cm}{	\rotatebox[origin=c]{90}{ alternating }		}
		&
		\tikzsetnextfilename{eFt_evol_alt} 
		\begin{tikzpicture}
			\begin{axis}[
						ymin=\ymin,		ymax=\ymax,
						ylabel={$\Delta F_t / \si{\newton}$}, 
						xmin=\xmin,		xmax=\xmax,
						xtick=\xTick,	xlabel={samples$^\prime$ / --}, 
						width=\SzW,		height=\SzH,
						xminorgrids,	minor x tick num=\NumMinorXtick,
						xmajorgrids,
						]
				\SigLbUbFromTableB{eF_RLS_alt.csv}{\ColorA}
				\SigLbUbFromTableB{eF_EnKF_alt.csv}{\ColorB}
				\SigLbUbFromTableB{eF_xEnKF_alt.csv}{\ColorC}

			\end{axis}
		\end{tikzpicture}\\
	\end{tabular}
	\caption{ Evolution of the error in the tangential force $\Delta F_t$  with the \KienzleMdl{} 	}
	\label{fig:eFt_evol}
\end{figure}
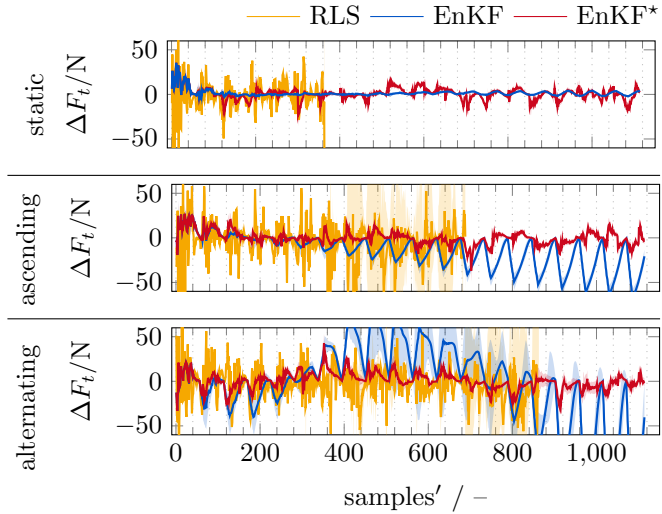

\tikzexternaldisable
\setlength{\tabcolsep}{6pt} 

Examining the evolution of the identified coefficients, \fig{ktmt_evol}, suggested an even worse performance of the \ac{RLS} than by just considering the model error. It was not capable of following the trending coefficients. 
However, the \ac{EnKF} seemed to follow the trend in the exponential coefficient $m_t$ quite well --  perfectly in the case of the linear ascent. This might have been due to the comparatively small change, since the exponential model coefficient has a narrow range and a \SI{20}{\percent} linear change did not require to maintain a large spread of the ensemble. However, it fails to follow the linearly ascending coefficient $k_t$ and revealed an extraordinary variance in the case of an alternating coefficient $k_t$ (light blue area). It appears that the changes are too drastic and too fast; in particular the alternating case of $k_t$ suggests a time-delay of the identification.
The repeated inflation of the ensemble (abbreviated here as \ac{EnKF}$^\star$, red lines) shows a small variance tube. The filter follows the trending coefficients and still exhibits good convergence in the static case. In general, the \ac{EnKF}$^\star$ leads to an unsteady convergence, which becomes especially obvious in comparison to the classic \ac{EnKF} in the static case.

\def\SzW{4cm}
\def\SzH{3.5cm}

\def\xmin{-10}
\def\xmax{1150}
\def\yminKi{1200}
\def\ymaxKi{2600}
\def\yminMi{0.1}
\def\ymaxMi{.42}
\def\xTick{0,500,1000,1500}
\def\NumMinorXtick{10}
\def\ColorA{wzl-yellow}
\def\ColorB{wzl-blue}
\def\ColorC{wzl-red}
\def\ColorD{black}

\setlength{\tabcolsep}{0pt} 
\tikzexternalenable

\begin{figure}[htp]
	\begin{tabular}{ l c c } 
		& \multicolumn{2}{r}{%
			\sampleline{\ColorA} \ac{RLS} \sampleline{\ColorB} \ac{EnKF} \sampleline{\ColorC} \ac{EnKF}$^\star$ \sampleline{\ColorD,style=dashed,thick} original
			}
		\\
		\raisebox{1cm}{	\rotatebox[origin=c]{90}{ static }		}
		&
		\tikzsetnextfilename{ktmt_evol_kt_sta} 
		\begin{tikzpicture}
			\begin{axis}[
						ymin=\yminKi,		ymax=\ymaxKi,
						ylabel={$\Delta k_t / \si{\newton\per\milli\meter}$}, 
						xmin=\xmin,		xmax=\xmax,
						xtick=\xTick,	xticklabels=\empty,
						width=\SzW,		height=\SzH,
						xminorgrids,	minor x tick num=\NumMinorXtick,
						xmajorgrids,
						]
				\SigLbUbFromTableB{ki_RLS_sta.csv}{\ColorA}
				\SigLbUbFromTableB{ki_xEnKF_sta.csv}{\ColorC}
				\SigLbUbFromTableB{ki_EnKF_sta.csv}{\ColorB}
				\SigFromTableSingle{ki_ogl.csv}{sta}{\ColorD,style=dashed,thick}
			\end{axis}
		\end{tikzpicture}
		&
		\tikzsetnextfilename{ktmt_evol_mt_sta} 
		\begin{tikzpicture}
			\begin{axis}[
						ymin=\yminMi,		ymax=\ymaxMi,
						ylabel={$\Delta m_t / -$}, 
						xmin=\xmin,		xmax=\xmax,
						xtick=\xTick,	xticklabels=\empty,
						width=\SzW,		height=\SzH,
						xminorgrids,	minor x tick num=\NumMinorXtick,
						xmajorgrids,
						]
				\SigLbUbFromTableB{mi_RLS_sta.csv}{\ColorA}
				\SigLbUbFromTableB{mi_xEnKF_sta.csv}{\ColorC}
				\SigLbUbFromTableB{mi_EnKF_sta.csv}{\ColorB}
				\SigFromTableSingle{mi_ogl.csv}{sta}{\ColorD,style=dashed,thick}
			\end{axis}		
		\end{tikzpicture}
		\\
		\hline 
		
		\raisebox{1cm}{	\rotatebox[origin=c]{90}{ ascending }		}
		&
		\tikzsetnextfilename{ktmt_evol_kt_asc} 
		\begin{tikzpicture}
			\begin{axis}[
						ymin=\yminKi,		ymax=\ymaxKi,
						ylabel={$\Delta k_t / \si{\newton\per\milli\meter}$}, 
						xmin=\xmin,		xmax=\xmax,
						xtick=\xTick,	xticklabels=\empty,
						width=\SzW,		height=\SzH,
						xminorgrids,	minor x tick num=\NumMinorXtick,
						xmajorgrids,
						]
				\SigLbUbFromTableB{ki_RLS_asc.csv}{\ColorA}
				\SigLbUbFromTableB{ki_EnKF_asc.csv}{\ColorB}
				\SigLbUbFromTableB{ki_xEnKF_asc.csv}{\ColorC}
				\SigFromTableSingle{ki_ogl.csv}{asc}{\ColorD,style=dashed,thick}
			\end{axis}
		\end{tikzpicture}
		&
		\tikzsetnextfilename{ktmt_evol_mt_asc} 
		\begin{tikzpicture}
			\begin{axis}[
						ymin=\yminMi,		ymax=\ymaxMi,
						ylabel={$\Delta m_t / -$}, 
						xmin=\xmin,		xmax=\xmax,
						xtick=\xTick,	xticklabels=\empty,
						width=\SzW,		height=\SzH,
						xminorgrids,	minor x tick num=\NumMinorXtick,
						xmajorgrids,
						]
				\SigLbUbFromTableB{mi_RLS_asc.csv}{\ColorA}
				\SigLbUbFromTableB{mi_xEnKF_asc.csv}{\ColorC}
				\SigLbUbFromTableB{mi_EnKF_asc.csv}{\ColorB}
				\SigFromTableSingle{mi_ogl.csv}{asc}{\ColorD,style=dashed,thick}
			\end{axis}		
		\end{tikzpicture}
		\\
		\hline 
		
		\raisebox{1.6cm}{	\rotatebox[origin=c]{90}{ alternating }		}
		&
		\tikzsetnextfilename{ktmt_evol_kt_alt} 
		\begin{tikzpicture}
			\begin{axis}[
						ymin=\yminKi,		ymax=\ymaxKi,
						ylabel={$\Delta k_t / \si{\newton\per\milli\meter}$}, 
						xmin=\xmin,		xmax=\xmax,
						xtick=\xTick,	xlabel={samples$^\prime$ / --}, 
						width=\SzW,		height=\SzH,
						xminorgrids,	minor x tick num=\NumMinorXtick,
						xmajorgrids,
						]
				\SigLbUbFromTableB{ki_RLS_alt.csv}{\ColorA}
				\SigLbUbFromTableB{ki_EnKF_alt.csv}{\ColorB}
				\SigLbUbFromTableB{ki_xEnKF_alt.csv}{\ColorC}
				\SigFromTableSingle{ki_ogl.csv}{alt}{\ColorD,style=dashed,thick}
			\end{axis}
		\end{tikzpicture}
		&
		\tikzsetnextfilename{ktmt_evol_mt_alt} 
		\begin{tikzpicture}
			\begin{axis}[
						ymin=\yminMi,		ymax=\ymaxMi,
						ylabel={$\Delta m_t / -$}, 
						xmin=\xmin,		xmax=\xmax,
						xtick=\xTick,	xlabel={samples$^\prime$ / --}, 
						width=\SzW,		height=\SzH,
						xminorgrids,	minor x tick num=\NumMinorXtick,
						xmajorgrids,
						]
				\SigLbUbFromTableB{mi_RLS_alt.csv}{\ColorA}
				\SigLbUbFromTableB{mi_xEnKF_alt.csv}{\ColorC}
				\SigLbUbFromTableB{mi_EnKF_alt.csv}{\ColorB}
				\SigFromTableSingle{mi_ogl.csv}{alt}{\ColorD,style=dashed,thick}
			\end{axis}		
		\end{tikzpicture}		
		\\
	\end{tabular}
	\caption{ Evolution of the tangential coefficients of the \KienzleMdl{}  	}
	\label{fig:ktmt_evol}
\end{figure}
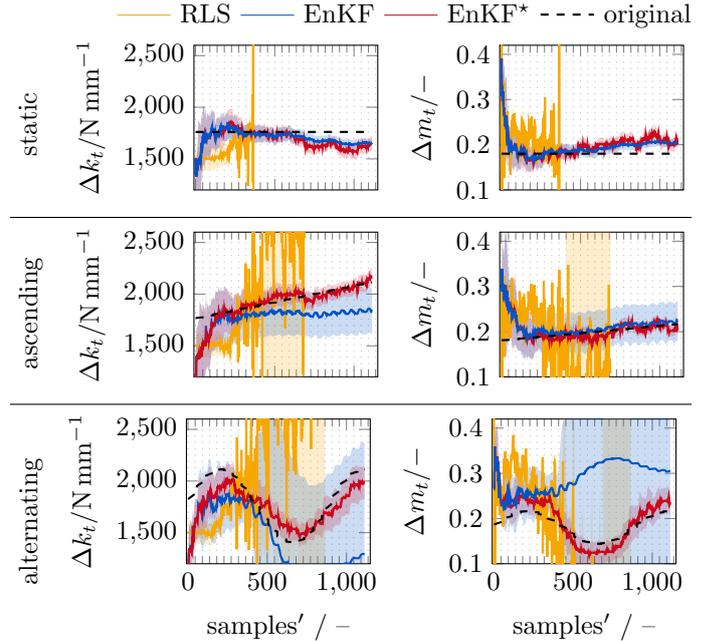

\tikzexternaldisable
\setlength{\tabcolsep}{6pt} 


\section{ Conclusion } 
A repeated inflation of (a subset of) the ensemble allows to identify time-varying coefficients. 
Instead of restarting the identification on a regular basis, the inflation maintains the ensemble mean as its current best-guess while enlarging the spread of the ensemble. This ensures that the ensemble explores the whole search-space at all times.
Inflating the ensemble consistently changes the subspace; therefore, convergence within the initial ensemble cannot be guaranteed anymore.
To regain control of the identification, it is important to impose box-constraints on the filter. Furthermore, the amount of inflation is set to a fraction of the variance of the initial ensemble to limit the spread. The fraction can be smaller the smaller the step-size of the ensemble inflation.
As an idea from the mean-field theory, only a subset of the ensemble is inflated in order to smooth the continuous identification.

The results show that the classic \ac{EnKF} is only able to follow trending coefficient in exceptional cases, \ie slow and small changes and with a high sensitivity to the initial ensemble. A repeated inflation of the ensemble drastically reduced this sensitivity and presented the only option to follow all cases of trending coefficients. However, at the cost of a slightly worse accuracy in the static case compared to the classic \ac{EnKF}. The \ac{RLS}, which served as a benchmark, revealed an immense proneness to the artificial white Gaussian noise in the measurements.

Future research will be placed on the integrating the mean field theory to the \ac{EnKF} in order to decrease oscillation and increasing stability through an adaptive time-step within the filter.
In milling, identification of a time-varying force model must be combined with the quasi-static identification of a model of the radial deviation of the tool. Eventually, the application to real measurement signals remain due.

\begin{ack}
	The authors would like to thank the German Research Foundation DFG for the kind support within the Cluster of Excellence \qml{}Internet of Production\qmr{} (Project ID: 390621612).
\end{ack}

\bibliography{ifacconf}

\end{document}